 \documentclass[11pt,reqno]{amsart}
\usepackage{amsmath,bm}
\usepackage{amsfonts}
\usepackage{amssymb}
\usepackage{amsthm}
\usepackage{amscd}
\usepackage{ytableau}
\usepackage{amsmath}
\usepackage{amssymb}
\usepackage{amsthm}
\usepackage{graphics}
\usepackage{eucal}
\usepackage{mathrsfs}
 \usepackage{amsmath,bm}
\usepackage{amsfonts}
\usepackage{amssymb}
\usepackage{amsthm}
\usepackage{amscd}
\usepackage{a4wide}
\usepackage{geometry}        
\usepackage{amsmath}
\usepackage{amssymb}
\usepackage{amsthm}
\usepackage{graphics}
\usepackage{eucal, times}
\usepackage{mathrsfs}
\usepackage{fancyhdr}
\usepackage[colorlinks=true, pdfstartview=FitV, linkcolor=blue, 
            citecolor=blue, urlcolor=blue]{hyperref}
\usepackage{memhfixc} 
\usepackage[utf8]{inputenc} % Any characters can be typed directly from the keyboard, eg ‚àö¬©‚àö√ü‚àö¬±
\usepackage{textcomp} % provide lots of new symbols
\usepackage{graphicx}  % Add graphics capabilities
\usepackage{flafter}  % Don't place floats before their definition
\usepackage{bm}  % Define \bm{} to use bold math fonts
\usepackage{tikz}
\usepackage{color}
\usepackage{memhfixc}  % remove conflict between the memoir class & hyperref
\usepackage{pdfsync}  % enable tex source and pdf output syncronicity
\usepackage{amssymb}
\usepackage{amsfonts}

\usepackage{bm}

\usepackage{tikz}
\usetikzlibrary{shapes.geometric} % for shape=ellipse
\usetikzlibrary{calc}

\usetikzlibrary{calc, positioning,fit,intersections}

\theoremstyle{plain}
\newtheorem*{maintheorem*}{Main Theorem}
\newtheorem*{thm*}{Theorem}
\newtheorem*{thma*}{Theorem A}
\newtheorem*{thmaa*}{Theorem A'}
\newtheorem*{thmb*}{Theorem B}
\newtheorem*{thmo*}{Theorem 1.1}
\newtheorem*{thmc*}{Theorem C}
\newtheorem*{thmd*}{Theorem D}
\newtheorem*{thmf*}{Theorem 4.1}
\newtheorem*{remark*}{Remark}
\newtheorem*{conjecture*}{Conjecture}
\newtheorem*{prop*}{Proposition}
\newtheorem*{lem*}{Basic Lemma}
\newtheorem{thm}{Theorem}[section]
\newtheorem{cor}[thm]{Corollary}
\newtheorem{lem}[thm]{Lemma}
\newtheorem{prop}[thm]{Proposition}

\newtheorem*{proofc*}{Proof of Theorem C}

\def\bbz{\mathbb{Z}}

\def\e{\mathbf{e}}

\def\x{\mathbf{x}}

\def\q{\mathbf{q}}

\usepackage{pgfplots}
\usepackage{pgfplots}
\pgfplotsset{width=7cm,compat=1.3}
\begin{document}
%%

%\thanks {\textit{Keywords: Diophantine approximation}}
\author{ Youssef Lazar }
\email{ylazar77@gmail.com}
\date{}

\title{Explicit solutions to the Oppenheim conjecture for  indefinite ternary diagonal forms}
\maketitle

\address{}

\begin{abstract}  We give a new proof of the Oppenheim conjecture for indefinite ternary diagonal forms of the type $x^{2}+y^{2} -\alpha  z^{2}$ where $ \alpha $ is an irrational number. Our method is explicit in the sense that we are able to construct a solution to the problem and we obtain an effective bound on the solution. The method is geometrical and is based on continued fractions. 
%% Text of abstract
  %% Text of abstract

\end{abstract}

%%%%%%%%%%%%%%%%%%%%%%%%%%%%%%%%%%%%%%%%%%%%%%%%%%%%%%%%%%%%%%%%%%%%%%   Introduction   %%%%%%%%%%%%%%%%%%%%%%%%%%%%%%%%%%%%%%%%%%%%%%%%%%%%%%%%%%%%%%%%%%%%%%%%%%%
\section{Introduction} \label{intro}

We are interested in the following diophantine problem, given any real number $\varepsilon >0$ and a positive irrational number $\alpha$, is there exists a nonzero vector $(x,y,z) \in \mathbb{Z}^{3}$ such that 
$$  \vert x^{2}+y^{2} -\alpha  z^{2}\vert \leq \varepsilon.$$

This apparently simple question found a solution only in the mid-eighties thanks to G.A. Margulis which proved that the answer is positive. In fact, Margulis showed in \cite{M1} a much more general statement which encompasses all indefinite quadratic forms in $n\geq 3$ variables, provided they are not proportional to a rational one.  This result was conjectured by Oppenheim in 1929 \cite{op1} and remained open in full generality until Margulis'  breakthrough. The Oppenheim conjecture reduces to the three dimensional case which is strangely the most difficult case. 
The strategy of the  proof used by Margulis, was to solve a particular case of another conjecture due to M.S. Raghunathan.  The resolution of the Oppenheim conjecture is a consequence of Margulis' proof of the Raghunathan conjecture in the case  $n=3$. Few time later, Ratner's proved the Raghunathan conjecture in full generality for any connected Lie group \cite{R1}. These results were the starting point of a tremendous amount of activity around which is now called \textit{homogeneous dynamics}. This point of view shows to be very fruitful in order to treat various unsolved problems especially in diophantine approximation.  The litterature about this conjecture and others related questions is abundent. The  interested reader may find most of the main contributions on this conjecture in Margulis' survey  \cite{M2} which is by far the most complete.   \\
A natural question is whether the Oppenheim conjecture could be proved with another method, namely without using homogeneous dynamics.  As far as we know the answer is negative for $n=3$, unless for a very specific case due to Watson which we will discuss later on. The most powerful method to solving diophantine inequalities is the \textit{Circle method} but it requires a large number of variables compared to the degree of the polynomial involved.  In the early times of the conjecture, Davenport and Heilbronn succeeded to prove the Oppenheim conjecture for irrational diagonal forms in five variables by using a variant of the Circle method \cite{DH}. Their proof has the advantage to be effective. The same result was proved earlier by Chowla for $ n \geq 9 $ using lattice points counting in irrational ellipsoids. 
 The barrier $ n=4$ has been breached by Oppenheim itself in its seminal paper \cite{op1} using some old results of Korkine and Zolotareff on representation of definite forms \cite{kz}. In the late seventies, Iwaniec \cite{Iwaniec} proved the Oppenheim conjecture for some quaternary diagonal forms using \textit{sieve theory}. A last attempt to prove the conjecture was due to R.K. Baker and H.P. Schlickewei who proved the conjecture in full generality for $ n \geq 21$ \cite{bk}. For quadratic forms, i.e. in degree 2, it seems that the circle method can only operate if $n \geq 5$. Using the full power of analytic methods combined with geometry of numbers,   an effective version of the Oppenheim conjecture was proved very recently for $n \geq 5$  by P. Buterus,  F. G\"{o}tze, T. Hille  and G.A~Margulis  \cite{M3}. This results have been sharpened by P. Buterus,  F. G\"{o}tze, T. Hille  in \cite{BGH} for diagonal forms extending Birch-Davenport method to dimensions at least five combined with a result of Schlickewei. The latter proofs are quite involved and very technical.

\vspace{0.5cm}

\noindent \textit{The three dimensional case.}\\
 It is noteworthy to mention the difficulty of the problem for $n=3 $. The case of forms $Q_{\alpha}(x,y,z)= x^{2}+y^{2} -\alpha  z^{2} $ we are concerned with shows a curious behaviour. Indeed it has been remarked by Eskin Margulis and Mozes (\cite{emm}, Theorem 2.2.) that $Q_{\alpha}(\bbz^{3})$ fails to be equidistributed for a dense set of values of $\alpha$. This contrasts with the analog in higher dimension, in the same paper it is proved that the set $Q(\bbz^{3})$ is equistributed in the real line given any form $Q$ of signature $(p,q)\neq (2,1)$ or $(2,2)$ which satisfy the assumptions of the Oppenheim conjecture.\\
For a very specific class of quadratics forms, Watson \cite{watson} gave an explicit proof of the Oppenheim conjecture by showing how to construct the solution and therefore providing bounds for the solution.

Watson considered quadratic forms of the type $Q(x,y,z)= x^{2}-a \alpha y^{2} - \alpha^{2}z^{2} $ where $a$ is a positive integer and $\alpha$ is an irrational number with continued fraction representation $[a; a, \ldots ]=[a; \overline{a}]$. When $a \geqslant 2$ such numbers are sometimes called \textit{silver means}, in analogy with the case $a=1$ which is just the golden ratio. The convergents of such numbers satisfies very a simple reccurence relation, if $c_{n}=p_{n}/q_{n}$ is the $n^{th}$ convergent of $\alpha$ then $q_{n}=p_{n-1}$. For each integer $n>0$, let us set 
\begin{center}
$  x_{n}=q_{n+1}$, $  y_{n}=q_{n}$ and $  z_{n}=q_{n-1}$.
\end{center}

\noindent  By means of easy manipulations Watson showed that
$$ \left| x_{n}^{2}-a\alpha y_{n}^{2} - \alpha^{2}z_{n}^{2}   \right| \leq \dfrac{\alpha + \overline{\alpha}}{q_{n} q_{n-1} B_{n} B_{n-1}} $$

\noindent where $\overline{\alpha}$ is the algebraic conjugate of $ \alpha $ and  $B_{n}=|\overline{\alpha} - p_{n}/q_{n}|$. Since $\alpha$ has bounded partial quotients, in fact all equal to $a$, the $B_{n}$'s are bounded. Thus, 
$$ \left| x_{n}^{2}-a\alpha y_{n}^{2} - \alpha^{2}z_{n}^{2}   \right| \ll_{n} \dfrac{1}{q_{n} q_{n-1}}. $$
Let us choose an arbitrary $\varepsilon >0$, then if $n$ is taken large enough in order to fullfill the inequality
$$  \dfrac{1}{q_{n} q_{n-1}} \leq   \dfrac{1}{q_{n-1}^{2}} \leq \varepsilon.$$
This ensures that $v_{n}= (q_{n+1}, q_{n}, q_{n-1})$ solves the Oppenheim conjecture for $n$ as above. Note that this gives an asymptotic sequence of solutions not only one solution.\\
A bound for the solution $v_{n}$ depends on the least integer $n_{1}$ such that $z_{n_{1}}=q_{n_{1}-1}=\dfrac{1}{\sqrt{\varepsilon}}$. Thus for $n\geq n_{1}$ 

\begin{equation}
\| v_{n} \|_{\infty}= q_{n+1} \ll \dfrac{1}{\sqrt{\varepsilon}}.
\end{equation}
\noindent This result is quite exceptional among the bunch of results surrounding the Oppenheim conjecture. In fact, it gives a computable solution and it is \textit{effective} in the sense that it gives a bound of the sequence $F(N)= 
\displaystyle \min_{v\in \bbz^{3}, v \neq 0, \| v\|_{\infty} < N} |Q(v)|$.  As we have seen, by taking $N=\varepsilon^{-1/2}$, Watson's result gives $$F(N)= 
\displaystyle \min_{v\in \bbz^{3}, v \neq 0, \| v\|_{\infty} < N} |Q(v)| \ll N^{-2}.$$

The problem of effectiveness in Margulis' theorem amounts to finding optimal bounds for $ F(N) $. Although, Ratner's theorems are not effective in general, Lindenstrauss and Margulis \cite{el} succeeded to overcome this issue by giving upper bounds on $F(1/\varepsilon)$ of the form $e^{P(1/\varepsilon)} $ for some polynomial $P$.
 Their deep result is valid for \textit{all} indefinite forms in degree three and is based on homogeneous dynamics.  Shorty after Bourgain \cite{bou} gave optimal bounds for $F(N)$ for ternary diagonal forms. The works of Ghosh- Gorodnik-Nevo \cite{GGN}, Ghosh-Kelmer  \cite{GK} and Athreya-Margulis \cite{AM} gave closely related results for generic families of quadratic forms. \\
  The bound provided in Watson's result is outstanding, in the sense that, as far as we know, this is the best known bound for an individual quadratic form. Indeed, one of the output of Bourgain's result predicts that, under Lindel\"{o}f hypothesis for the Riemann Zeta function the best bound one can hope for a generic form is $F(N) \ll N^{-1+o(1)}$. Watson's peculiar example improves it by a factor $N^{-1}$. Note a slight difference with Bourgain, indeed he considered forms of the type $ Q(x,y,z)=x^{2}_{1} + \alpha_{2} x_{2}^{2} - \alpha_{3}x_{3}^{2} $ with $\alpha_{2},\alpha_{3}>0$ whereas Watson's example is of the form $ Q(x,y,z)=x^{2}_{1} - \alpha_{2} x_{2}^{2} - \alpha_{3}x_{3}^{2} $ with $\alpha_{2},\alpha_{3}>0$. Be that as it may, Watson's result is the best result one can expect in solving diophantine inequalities of the form $|Q(v)| \leq \varepsilon$. 

\begin{center}
\begin{figure}

\begin{tabular}{|c|c|c|c|c|c|}
\hline

Method & Variables(s)&Type &Quantative & Effective    & Explicit      \\\hline

Homogeneous  Dynamics & $n \geq 3$ &general &$\checkmark $ & $\checkmark $    &$ \times$  \\\hline

Circle method & $n \geq 5$ & general & $\checkmark $ & $\checkmark $    &$ \times$  \\\hline

Geometry of  numbers &$n=4$, $n \geq 9$ &diagonal &$\checkmark $ & $\checkmark $    &$ \times$  \\\hline

Sieve Theory & $n =4$ & diagonal &$\times $ & $\times $    &$ \times$  \\\hline

Continued Fractions & $n =3,4$ &diagonal &$\times $ & $\checkmark $    &$ \checkmark$  \\\hline

\end{tabular}
\caption{Comparison of the different proofs of the Oppenheim conjecture.}
\end{figure}
\end{center}

\subsubsection*{The main results}
The aim of the paper is to construct explicit solutions to the Oppenheim conjecture for ternary forms of the type $Q_{\alpha}(x,y,z)=   x^{2}+y^{2} -\alpha  z^{2}$ where $ \alpha \notin \mathbb{Q} $.  In turn, one is able to obtain effective bounds on such solutions. The proof essentially relies on  diophantine properties of the irrational number  $\beta = \sqrt{\alpha}$, more precisely its measure of irrationality.  The measure of irrationality of an real number $\beta$ is defined as the least positive real number $\mu$ such that for 
\begin{center}
$\mu(\beta) = \inf  \{ \omega \ : \  \dfrac{C}{q^{\omega +\sigma}} < \left| \beta - \dfrac{p}{q} \right| $ for all rational $p/q$ ($q >0$), every $\sigma >0$ and for some constant $C>0 \}$.
\end{center} 

\noindent A deep theorem due to Roth states that $\mu(\beta)=2$ whenever $\beta$ is an algebraic number. The converse is not true, indeed $\theta(e)=2$ whereas the constant $e \approx 2.718$ is a transcendental number.  There exists transcendental numbers $x$ for which $\mu(x)=\infty$, these are termed \textit{Liouville numbers}.

\noindent Let us fix an arbitrary small parameter $\sigma >0$. 
For any irrational number $\beta$ which is not a Liouville number, we define the following quantity 
$$ \theta_{\sigma}(\beta):= \mu(\beta)-1 + \sigma. $$

\noindent In any case, one has $\theta > 1$ and the following diophantine condition holds for any irrational $\beta$

\begin{equation}
\inf_{q \geq 1} q^{\theta} \langle q \beta \rangle > 0
\end{equation}
           
\noindent where $\langle x \rangle$ denotes the distance of a real number $x$ to the nearest integer. We are going to give an explicit proof of the  Oppenheim conjecture for quadratics forms  of the type $Q_{\alpha}(x,y,z)= x^2 + y^2 - \alpha z^2$ where $\alpha$ is an irrational number. The convergents of $\beta$ are simply denoted $\textbf{c}_{n}= \textbf{p}_{n}/ \q_{n}$ and $\theta$ stands for $\theta_{\sigma}(\beta)$.
\begin{thm} \label{main} Given any real number $\varepsilon >0$ and a positive irrational number $\alpha$. There exists a nonzero vector $v=(x,y,z) \in \mathbb{Z}^{3}$ such that 
$$\vert  Q_{\alpha}(v) \vert  \leq \varepsilon.$$
Moreover, if $\beta$ is not a Liouville number then the solution satisfies
$$  \| v \|_{\infty} \ll  \q_{2n_{1}}^{2/(\theta + 1)}$$

\noindent where $\q_{2n_{1}}$ is the denomimator of the convergent  of order $2n_{1}$ of $\beta$ with
$$ n_{1}(\varepsilon) = 2+  \lfloor \dfrac{\theta +1}{\theta-1}  \left| \ln \left( \varepsilon \right) \right| /\ln 2 \rfloor.$$

\end{thm}

We can extend the class of forms for which the conjecture is valid by considering classes of forms equivalent to the type $Q_{\alpha}$ as  given in Theorem \ref{main}.  
Given a subgroup $N$ of $\mathrm{GL}(3,\mathbb{R})$, we say that two quadratic forms $Q_{1}$  and $Q_{2}$ are $N$-equivalent if there exists a $g  \in N$ such that $Q_{1}(x)= Q_{2}(gx)$ for every $x \in \mathbb{R}^{3}$. Any indefinite ternary form $Q$ is $\mathrm{SL}(3,\mathbb{R})$-equivalent to $Q_{0}(x,y,z)= x^{2} + y^{2} -z^{2}$. Let us denote by $H$ the subgroup of $\mathrm{GL}(3,\mathbb{R})$ defined by
$$H =\left\lbrace  \left[\begin{array}{c|c}A & 0 \\\hline 0 & h_{33}\end{array}\right] :  A \in \mathrm{SL}(3,\mathbb{Q}) , h_{33} \notin \mathbb{Q} \right\rbrace.$$
From Theorem \ref{main} we derive the following result.
\begin{cor} \label{cor}
\begin{enumerate}
\item Suppose that $Q$ is an indefinite quadratic form which is $\mathrm{SL}(3,\mathbb{\mathbb{Q}})$-equivalent to a $Q_{\alpha}$ with $\alpha \notin \mathbb{Q}$. Then the Oppenheim conjecture holds for $Q$. \\

\item Suppose that $Q$ is an indefinite quadratic form  which is $H$-equivalent to  $Q_{0}$. Then the Oppenheim conjecture holds for $Q$. In particular if $ Q_{1}=f(x,y)- \beta^{2} z^{2}$ where $f(x,y)$ is a rational binary form and $\beta \notin \mathbb{Q}$ then the conjecture holds for $Q_{1}$.

 %Moreover, if $Q=Q_{0}^{h}$ where $$h = \left[\begin{array}{c|c}A & 0 \\\hline 0 & h_{33}\end{array}\right] \in H. $$
%where $h_{33}$ is not a Liouville number, then for every $\varepsilon >0$ there exists an integral solution of $|Q(v)| \leq \varepsilon$ with norm
%$$   \| v \|_{\infty} \ll \| A^{-1}\| \q_{2n_{1}}^{2/(\theta + 1)}$$ 
%with $\theta$ and $n_{1}$ as in Theorem \ref{main} with $\beta=h_{33}$.
\end{enumerate}

\end{cor}

\noindent \textbf{Remarks.} (1) A great advantage of our method is that we know how to construct the solution.  As a byproduct we obtain an effective  bound on the solution. The quality of the bound depends on the value of $\theta$ and the growth of the denominators $(\q_{n}(\beta))_{n \geqslant 1}$ of the convergents of $\beta$.

(2) The idea of the proof is geometrical and relies on the following observations. The line of equation $x=\beta z$ is a generatrix for the cone $\{Q_{\alpha}=0 \}$ restricted to the plane $y=0$. For every $\varepsilon >0$, this line is inside the region $\{-\varepsilon \leq Q_{\alpha} \leq \varepsilon  \}$ and because $\beta$ is irrational, this line cannot contain a nontrivial lattice point. Nevertheless, Dirichlet's approximation theorem tells us that there exists lattice points lying arbitrarily near the line at any level of precision. Given any $\varepsilon >0$,  one expects that such lattice point lies in the region $\{|Q_{\alpha}| \leq \varepsilon \}$. Unsurprinsingly we show in section \S 2 that Dirichlet's theorem is not enough to prove the Oppenheim conjecture for $Q_{\alpha}$.  To overcome this problem we introduce a sequence of rational lines which are nearly parralel to the line passing through a lattice point $ u_{n}=(x_{n}, 0, z_{n}) $ given by Dirichlet's given a certain order of approximation $\q_{2n}^{-1}$, i.e.
$$  \mathrm{dist}(u_{n}, \mathcal{L}_{\beta}) \ll \dfrac{1}{\q_{2n}^{1-\eta}} $$

where  $ 1 \leq z_{n} \leq \q_{2n}^{1-\eta} $ with   $ \eta = \dfrac{\theta -1}{\theta+1}$. Given any $n$, we define the line $\mathcal{L}_{\beta}^{n}$  by setting
$$ (\mathcal{L}_{\beta}^{n})  : ~u_{n} + \mathbb{R}(c_{2n}(\beta), \dfrac{1}{\q_{2n}(\beta)} ,1). $$
 \\
A parametrization of this line for the downward direction is given by 
$$(\mathcal{L}_{\beta}^{n})^{+} ~:~ v_{n}(t) = (x_{n}-t c_{2n}(\beta), -\dfrac{t}{\q_{2n}(\beta)}, z_{n}-t) ~~( t \geq 0) .$$ 

The proposition \ref{exit times} is going to show that for $n$ large enough the parametrization $v_{n}(t)$ of the intersection ${\mathcal{L}_{\beta}^{n}}^{+} \cap \{|Q_{\alpha}| \leq \varepsilon \}$ is supported by two disjoints intevals $I_{n}^{1}$ and $I_{n}^{2}$. Thus, in order to have a lattice point in $ \{|Q_{\alpha}| \leq \varepsilon \}$ it suffices to find a multiple of $\q_{2n}$, say $t_{n}$, in the union of $I_{n}^{1}$ and $I_{n}^{2}$. In this case, one can clear denominators and the solution is given by $v_{n}(t_{n}) \in \bbz^{3} \cap {\mathcal{L}_{\beta}^{n}}^{+} \cap \{|Q_{\alpha}| \leq \varepsilon \}$. The key lemma \ref{multiple} says that this is possible if $n$ is greater or equal than some integer $n_{1}(\varepsilon)$ which is explicitely computable.\\
3) When $\beta=\sqrt{\alpha}$ is a Liouville number, we can easily prove that the Oppenheim conjecture is satisfied in dimension $n=2$ for the form $q(x,z)=x^{2}-\beta^{2}z^{2} $.  Since we have $Q_{\alpha}(x,0,z)=q(x,z)$, then the Oppenheim conjecture is satisfied for $Q_{\alpha}$. \\
4) In the case when $\beta$ is not a Liouville number, it is  always possible to find a real number $\theta>1$ large enough such that for every integer $q \geqslant 1$, 
$$ q^{\theta} \langle q \beta \rangle > 0.$$
The irrationality measure  $\mu(\beta)$ is introduced only with the aim of obtaining optimal bounds and to quantify the growth of the sequence $\q_{n+1}(\beta)/ \q_{n}(\beta)$. For instance the main theorem gives an explicit integral solution for our favorite example of irrational indefinite form $Q(x,y,z)= x^{2}+ y^{2}- \sqrt{2} z^{2}$.

5) The corollary \ref{cor} shows that we can find a solution to the Oppenheim problem for quadratic forms of the type
$$Q(x,y,z)= ax^{2}+bxy+ cy^{2}- \alpha z^{2}$$
where $a,b,c \in \mathbb{Q}$ and $\alpha \notin \mathbb{Q}$. This is the best we can do, and it would be interesting to find explicit solutions for general indefinite irrational forms. For the general case, one would be led to use the mutidimensional version of the Dirichlet's approximation theorem. Using the same kind of strategy applied to a product of linear forms instead of a quadratic form, the author was able to derive a set of sufficent conditions for the Littlewood conjecture to hold.

\section{Sequences of Rational Lines intersecting $\{|Q_{\alpha}| \leq \varepsilon \}$ } 
We focus our attention on forms of the type 
$$Q_{\alpha}(x,y,z) = x^{2}+ y^{2} -\alpha z^{2} $$
where $\alpha \in  \mathbb{R}_{+}$. We assume that $\alpha$ is irrational and therefore the form $Q_{\alpha}$ is an indefnite quadratic form which is not proportional to a form which rational coefficients.  The output of Margulis's result tells us that for every $\varepsilon >0$, there must exist a nonzero lattice vector $v \in \mathbb{Z}^{3}$ such that 
\begin{equation}\label{opp}
0\leq |Q_{\alpha}(v)| \leq \varepsilon.
\end{equation} 

We are going to reprove this result by constructing an explicit solution to this problem, i.e. to find a nonzero vector in $\mathcal{A}(\varepsilon) \cap \mathbb{Z}^{3}$ where he domain $\mathcal{A}(\varepsilon)$ is delimited by the level sets 
$\{ Q_{\alpha}=-\varepsilon \}$ and $\{ Q=\varepsilon \}$.\\
 A parametization of the cone $\{ Q_{\alpha}=0 \}$ is as follows, 
\begin{equation} \left\lbrace \begin{array}{cc}
x(t, \theta)&=   \sqrt{\alpha} ~t \cos \theta \\
 y(t, \theta)&=   \sqrt{\alpha} ~t \sin \theta\\
 z(t,\theta)&=   t. \ \ \ \ \ \  \ \ \ \ 
\end{array}\right.  (0 \leqslant \theta < 2\pi).
\end{equation}
This parametrization shows that the cone $\{ Q_{\alpha} =0 \}$ is generated by a continous family of lines given by $\mathcal{L}_{\alpha}(\theta) = \mathbb{R}(\sqrt{\alpha}~ \cos \theta,\sqrt{\alpha}~ \sin \theta, 1 )$ where the angle $\theta$ varies in $[0,2\pi)$. The line corresponding to the intersection of the $xz$-plane with the cone $\{ Q_{\alpha} =0 \}$ is exactly given by $\mathcal{L}_{\alpha}(0) = \mathbb{R}(\sqrt{\alpha},0, 1 )$, we denote it by $\mathcal{L}_{\beta}$ where $\beta=\sqrt{\alpha}$. An equation of this line in the $xz$-plane is just $x= \beta z$. Since $\beta^{2}=\alpha$ is irrational, $\beta$ itself is irrational too. Thus given any positive integer $N>1$ we obtain from  Dirichlet's Theorem
that there exists $(p_{0},q_{0}) \in \mathbb{N}^{2}$ with $1 \leqslant q_{0} \leqslant N$ such that 
\begin{equation}\label{dirichlet}
 |p_{0}-\beta q_{0}| \leqslant \frac{1}{N}.
\end{equation}
This tells us that we can always find a lattice vector $(p_{0}, 0, q_{0})$ arbitrarily near the line of equation $x= \beta z$ in the $xy$-plane provided $\beta$ is irrational.

\subsection{ Irrationality Measures} We follow the notations of \cite{hata}, section 3. \\
For each real number, let $\langle x\rangle$ denote the distance of $x$ to the closest integer. Dirichlet's theorem says that $ \inf_{q \geq 1} q\langle q\beta\rangle < 1$, the question is to know in which extend one can improve this approximation. We can assign to $\beta$ a number called the irrationality measure of $\beta$ which is defined as follows,
\begin{center}
$\mu(\beta) :=  \inf  \{ \omega \in \mathbb{R}_{+}   \ : \  \inf_{q \geq 1} q^{\omega -1 + \sigma}  \langle q \beta\rangle >0 $ for every  real $\sigma >0 \}$.
\end{center} 

In other words, 
\begin{center}
$\mu(\beta) = \inf  \{ \omega \ : \  \dfrac{C}{q^{\omega +\sigma}} < \left| \beta - \dfrac{p}{q} \right| $ for all rational $p/q$ ($q >0$), every $\sigma >0$ and for some constant $C>0 \}$
\end{center}

\textit{An alternative definition of $\mu$}\\
 Suppose $\beta$ has an infinite continued fraction expansion $\beta=[b_{0}; b_{1}, b_{2}, \ldots  ]$, the $n^{th}$ convergent of $\beta$ is the rational number $\mathbf{c}_{n}(\beta)=[b_{0}; b_{1}, \ldots, b_{n}  ]$ which has  reduced expression $\dfrac{\mathbf{p}_{n}(\beta)}{\q_{n}(\beta)}$.  	
Then the measure of irrationality of $\beta$ is related to the growth of the denominators of $\mathbf{c}_{n}(\beta)$ through the following relation which can be taken as an alternative defintion of $\mu$,
$$ \mu(\beta) = 1 + \limsup_{n} \dfrac{\ln \q_{n+1}(\beta)}{\ln \q_{n}(\beta)}.$$

\noindent Provided the existence of the limit, one has the following asymptotic behaviour
$$ q_{n+1} \asymp q_{n}^{\mu -1}. $$

\noindent If we denote by $ \lambda_{n}$ the ratio $ q_{n+1}/q_{n}$, the last asymptotic estimate could be read as follows
\begin{equation}\label{lambdamu}
\lambda_{n} \asymp q_{n}^{\mu -2}.
\end{equation}

Its lowest value for an irrational number is $\mu(\beta)=2$ and it is reached for any algebraic number of degree $d \geq 2$. This fact is a highly non trivial theorem due to Roth \cite{roth}.  In the other extreme side, the value $\mu(\beta)=\infty$ correspond to the case when  $\beta$ is the Liouville number. In general, it is extremely difficult to compute this measure in practice. \\
A nice consequence of Roth's theorem is that $x$ is a transcendental number as soon as $\mu(x)>2$.  Unfortunately, this criterion is not enough to characterize  transcendental numbers because  the converse of Roth's result is not true. Indeed, Adams' proved that $\theta(e)=2$ showing that a transcendental number could reach the same bound (see e.g. \cite{davis}). More precisely, it can be proved that for all rational numbers $p/q$ ($q \geq 2$)
$$ | qe -p  |> c_{1} \dfrac{\log \log q}{q\log q} . $$
Since the continued fraction expansion of $e$ is given $e=[2;1,2, \overline{1,1,2n}]^{n\geqslant 2}$,  its partial quotients are unbounded. This implies that $e$ is not a badly approximable number, thus $\inf_{q \geq 1} q\langle qe\rangle =0$.  But for every $\sigma>0$, it is not difficult to see that for every $q \geqslant 2$
$$ \dfrac{\log \log q}{\log q} > \dfrac{1}{q^{\sigma}}.$$
This shows that for every $\sigma>0$ and rationals $p/q$
$$| qe -p| > \dfrac{c}{q^{1+\sigma}}$$
for some constant $c$. The latter amounts to say that $\mu(e)=2$, and it shows that some transcendental numbers are not well-approximated by rationals and behave like algebraic numbers in view of Roth's theorem.

\textit{The exponent theta associated to $\beta$}\\
By definition suppose that $\beta$ is not a Liouville number i.e. $\mu(\beta)< \infty$. Then for every $\sigma >0$ there exists $C>0$ such that for any $p,q$ integers with $q \geq 1$, 

\begin{equation}\label{diri}
 \dfrac{C}{{q}^{\mu(\beta) + \sigma}} < \left| \beta - \dfrac{p}{q}  \right|
\end{equation}

\noindent or also, 

\begin{equation}\label{diri}
 \dfrac{C}{{q}^{\mu(\beta) -1 + \sigma}} < \left| q\beta - p  \right|.
\end{equation}
\noindent Let us fix a real parameter $\sigma >0$ and introduce the following useful quantity associated with any irrational number $\beta$
$$ \theta(\beta) := \mu(\beta)+ \sigma-1.$$

\noindent This exponent gives a lower bound for the approximation of the irrational number $\beta$ by rational numbers provided it is not a Liouville number. It particular, since $\beta$ is not a Liouville number one has that  $C= \inf_{q \geqslant 1} q^{\theta}\|q \beta \|$ is positive and therefore for every $p$ and $q$ integers, $q \geq 1$

\begin{equation}\label{diri}
 \dfrac{C}{{q}^{\theta+1}} < \left| \beta - \dfrac{p}{q}  \right|.
\end{equation}

\subsection{Dirichlet versus Oppenheim}\noindent Dirichlet's theorem does not give a very precise estimate about how close is the lattice point $u_{0}=(p,0,q)$, obtained in (\ref{dirichlet}), to the line $\mathbb{R}(\beta, 0, 1)$. 
Let us explain why $u_{0}=(p_{0},0,q_{0})$ falls out $\mathcal{A}(\varepsilon)$ for any choice of $N$.  Otherwise the conjecture would be proved for $Q_{\alpha}$ and $u_{0}$ would be the solution. It is not difficult to quantify by how much Dirichlet's fails to prove the Oppenheim conjecture for $Q_{\alpha}$.

\noindent
In particular, combining (\ref{diri}) with Dirichlet approximation (\ref{dirichlet})  we have
$$ \dfrac{C}{q_{0}^{\theta+1}} < \left| \beta - \dfrac{p_{0}}{q_{0}}\right| \leqslant \dfrac{1}{q_{0}N}. $$

\noindent Taking the inverse if necessary, we can assume that  $u_0$ is in the first octant with $q_{0} \leqslant N$ one infers that
\begin{equation}\label{dirichletpositive}
 \dfrac{C}{q_{0}^{\theta}} < p_0 - \beta q_{0}   \leqslant 
 \dfrac{1}{N}.
\end{equation}

\noindent and the latter inequality gives in addition sharp bounds for $q_{0}$
\begin{equation}\label{q0}
(CN)^{1/\theta} < q_{0} \leqslant N.
\end{equation}

\noindent From (\ref{dirichletpositive})  we get

\begin{equation}\label{ineqqpbad}
\beta q_{0} + \dfrac{C}{q_{0}^{\theta}}< p_{0} < \beta q_{0} +\frac{1}{N}.
\end{equation}

\noindent Using (\ref{q0}) one obtains

\begin{equation}\label{ineqqpbad2}
 \beta (CN)^{1/\theta} + \dfrac{C}{N^{\theta}}< p_{0} \leq \beta N +\frac{1}{N}.
\end{equation}

\noindent Therefore 

\begin{equation}\label{ineqbadplus}
2 \beta (CN)^{1/\theta}+ \dfrac{C}{N^{\theta}}< p_{0} + \beta q_{0} \leq 2\beta N +\frac{1}{N}.
\end{equation}

\noindent We finally obtain the following bounds for $Q_{\alpha}(u_{0})=  p_{0}^{2}-  \beta^{2} q_{0}^{2}$

\begin{equation}\label{ineqbadfinal}
\frac{C}{N^{\theta}}\left(2 \beta (CN)^{1/\theta}+ \dfrac{C}{N^{\theta}} \right)< p_{0}^{2}-  \beta^{2} q_{0}^{2} \leq \frac{1}{N}\left(2 \beta N+ \dfrac{1}{N} \right).
\end{equation}

\noindent In particular we can do than the inequality $ Q_{\alpha}(u_{0}) < 2\beta + \dfrac{1}{N^{2}} $.   Thus Dirichlet's theorem is unable to provide a solution to the Oppenheim conjecture whatever the choice of  $N$.

\subsection{Error in the approximation by the convergents} We have a precise of the rate of error of this approximation, set $\e_{n}(\beta): = \beta - c_{n}(\beta)$, so we have (see e.g. Exercise 3.1.5. \cite{EW})
\begin{equation}\label{error1}
 \dfrac{1}{2\q_{n+1}(\beta)^{2}} \leq    |\e_{n}(\beta)| \leq  \dfrac{1}{\q_{n}(\beta)\q_{n+1}(\beta)} <\dfrac{1}{\q_{n}(\beta)^{2}}.
\end{equation}
The sequence $\q_{n}(\beta)$ is increasing and the rate of convergence is determined by the diophantine properties of $\beta$, in particular it tends to infinity with at least exponential rate since $2^{(n-2)/2} \leqslant \q_{n}(\beta)$ (\cite{Kh}, Theorem 12).  

%Since the partial quotients $(b_n)$ are bounded by a common value $M$, thus $\q_{n}(\beta) \leqslant (M+1)^{n}$. To sum up, $\beta$ beeing badly approximable has the fundamental property 
%\begin{equation}\label{boundqn}
%2^{n-1} \leqslant \q_{2n}(\beta) \leqslant (M+1)^{2n}
%\end{equation}

\noindent The convergents  $c_{n}(\beta)$ tends to $\beta$ by  oscillating so that the sign of $\e_{n}(\beta)$ is alternating. From now on, we choose even indices which implies that the error terms assume only positive values.  We infer that,

 \begin{equation}\label{error2}
\dfrac{1}{2\q_{2n+1}(\beta)^{2}} \leq    \e_{2n}(\beta) < \dfrac{1}{\q_{2n}(\beta)^{2}}.
\end{equation}

\subsection{Rational lines of approximation.} We introduce an object which is at the core of our strategy. It is a sequence of rational lines which will cross $\mathcal{A}(\varepsilon)$ in a sufficently large time in order to contain a lattice point. We have two degrees of freedom given by the integral parameters $n$ and $N$. We are going to reduce to merely one parameter, namely $n$. To do this, let us first fix the real parameter $$ \eta: = \dfrac{\theta -1}{\theta+1}$$  
where $$ \theta(\beta) = \mu(\beta)+ \sigma-1.$$
In all cases, $\theta >1$, and therefore
$$0< \eta<1.$$  
Let us choose $N$ to be a sequence $\left(N_{n} \right)$ satisfying the growth condition
\begin{equation}\label{growth}
  N_{n} =\q_{2n}^{1-\eta}.
\end{equation}

\noindent For each nonnegative integer $n$, Dirichlet's theorem tells us that there exists a two-dimensional lattice vector $(x_{n},z_{n})$ with $1 \leqslant z_{n} \leqslant N_{n}$ such that 
\begin{equation}\label{dirich1}
\left| \beta z_{n} - x_{n} \right| \leq \dfrac{1}{N_{n}}.
\end{equation}
Moreover, using  $\theta $ there exists a constant $C$ such that 
\begin{equation}\label{dirich2}
 \dfrac{C}{N_{n}^{\theta}} \leq \dfrac{C}{z_{n}^{\theta}} < \left| \beta z_{n} - x_{n} \right| \leq \dfrac{1}{N_{n}}.
\end{equation}

\noindent This gives the crucial bound on the denominators,
\begin{equation}\label{boundszn}
(CN_{n})^{1/\theta} <  z_{n} \leq  N_{n}
\end{equation}

 For each $n$, from (\ref{dirich1}) we form the three-dimensional integral vector $u_{n} := (x_{n}, 0, z_{n}) \in \mathbb{Z}^{3}$ which is close to the axis $x=\beta z$. As we have seen earlier Dirichlet's approximation theorem is not enough in order to ensure that $ u_{n} $  is in $\mathcal{A}(\varepsilon)$. However, we have at our disposal a sequence of lattice points $(u_{n})_{n}$ near $\mathcal{A}(\varepsilon)$ from which we built a sequence of affine lines $\mathcal{L}_{\beta}^{n}$  by setting
$$ (\mathcal{L}_{\beta}^{n})  : ~u_{n} + \mathbb{R}(c_{2n}(\beta), \dfrac{1}{\q_{2n}(\beta)} ,1).$$
\noindent  The lines $(\mathcal{L}_{\beta}^{n})$ are good candidates for containing lattice points in $ \mathcal{A}(\varepsilon)$.  Indeed the first interesting feature is that this lines pass through lattice points, namely the $u_n$'s, and such lines are directed by rational vectors so that they can contain latiice points.  Another crucial feature is geometrical, the fact that the lines are nearly parralel to the generatrix of the cone, namely the line $\mathbb{R}(\beta, 0,1)$. This+ leads us to expect that $(\mathcal{L}_{\beta}^{n})$ spends a sufficent amount of time in $ \mathcal{A}(\varepsilon)$ for $n$ large enough.

\noindent We will rather focus on the downward half-line  parametrized as follows $$(\mathcal{L}_{\beta}^{n})^{+} ~:~ v_{n}(t) = (x_{n}-t c_{2n}(\beta), -\dfrac{t}{\q_{2n}(\beta)}, z_{n}-t) ~~( t \geq 0) .$$ 
 We are interested in the intersection of this half-line with the domain $\mathcal{A}(\varepsilon)$. Note that for each increment of the index $n  $, the line will never remains in a same plane, in that two successive lines $(\mathcal{L}_{\beta}^{n})^{+} $ and $(\mathcal{L}_{\beta}^{n+1})^{+} $  will never be coplanar. A geometric observation allows us to guess that this line $(\mathcal{L}_{\beta}^{n})^{+}$ will cut the boundary of $\mathcal{A}(\varepsilon)$, namely $\{ Q_{\alpha} = \pm \varepsilon\}$, in at most four points. This will be made explicit in our computations. Our first task is to estimate the time spent by $(\mathcal{L}_{\beta}^{n})^{+}$ in $\mathcal{A}(\varepsilon)$. The answer is given in the following proposition.

\begin{prop} \label{exit times} Let $I^{n}$ be the set of times at which the half-line $(\mathcal{L}_{\beta}^{n})^{+}=\{ v_{n}(t)  \mid  t \geq 0)\}$ intersects $\mathcal{A}(\varepsilon)$. Then,  for $n$ large enough, $I$ is the union of two intervals  $I_{1}^{n}$ 
and  $I_{2}^{n}$.
\end{prop}

\noindent \textit{Proof.} Let us fix $\varepsilon >0$.  For each positive integer $n$,  the half-line $(\mathcal{L}_{\beta}^{n})^{+}$ lies in $\mathcal{A}(\varepsilon)$ if and only if for every $ t \geq 0 $
$$v_{n}(t) =(x_{n}-t c_{2n}(\beta), -\dfrac{t}{\q_{2n}(\beta)}, z_{n}-t) ~ \in \mathcal{A}(\varepsilon).$$ 
This amounts to say that the time variable $t$ is constrained to satisfy the inequalities
$$ -\varepsilon \leq (x_{n}-tc_{2n}(\beta))^{2} + (t \q_{2n}(\beta)^{-1})^{2} -\alpha (z_{n}-t)^{2} \leq \varepsilon.$$

\noindent Let us define the quadratic polynomial in the real variable $t$ (the time)
 $$f_{n}(t) :=  (x_{n}-tc_{2n}(\beta))^{2} + (t \q_{2n}(\beta)^{-1})^{2} -\alpha (z_{n}-t)^{2}.$$  

\noindent Ordering the terms we get
$$f_{n}(t)= \{c_{2n}(\beta)^{2}+ \q_{2n}(\beta)^{-2} - \beta^2  \}t^2 -2\{c_{2n}(\beta)x_{n}-y_{n}\beta^{2} \} t +\{  x_{n}^{2}-\beta^{2}z_{n}^{2}\}. $$

\noindent Which is important to us is the intersection points of the graph of $f_{n}(t)$ with the two lines corresponding to $\pm \varepsilon$. Thus we are reduced to solve the two following equations ( remember $\beta^2 = \alpha$ ) provided such solutions exists
$$f_{n}(t)  \pm \varepsilon = \{c_{2n}(\beta)^{2}+ \q_{2n}(\beta)^{-2} - \beta^2  \}t^2 -2\{c_{2n}(\beta)x_{n}-y_{n}\beta^{2} \} t +\{  x_{n}^{2}-\beta^{2}z_{n}^{2} \pm \varepsilon\} =0. $$

\noindent Set $A_{n} = c_{2n}(\beta)^{2}+  \q_{2n}(\beta)^{-2} - \beta^2 $, $B_{n} = -2\{c_{2n}(\beta)x_{n}-z_{n}\beta^{2} \} $ and $C_{n}^{\pm}=x_{n}^{2}-\beta^{2}z_{n}^{2} \pm \varepsilon $.\\

\noindent Thus one has to solve the (two) equations
$$A_{n} t^{2} + B_{n} t +C_{n}^{\pm}=0. $$
We need to estimate the discriminants $\Delta_{n}^{\pm}(\varepsilon)= B^{2}_{n}- 4A_{n}C_{n}^{\pm}$ and in fact we only need to focus on the roots and their relative distance not on their ordering nor their signs.

Since we have a nice control of the error in the approximation by the convergents, we replace $c_{2n}(\beta)$ by $\beta - \e_{2n}(\beta)$. The coefficients are therefore given by,

\begin{equation}\left\lbrace \begin{array}{ccc}
A_{n}&=  & -2\beta \e_{2n}(\beta)+  \e_{2n}(\beta)^{2} + \q_{2n}(\beta)^{-2}. \\
 B_n &=  &- 2(\beta x_{n} - x_{n}\e_{2n}(\beta) -z_{n}\beta^{2})= -2(\beta(x_{n}-\beta z_{n})-x_{n}\e_{2n}(\beta))\\
C_{n}^{\pm}&= &(x_{n}-\beta z_{n})(x_{n}+\beta z_{n}) \pm \varepsilon.
\end{array}\right.  
\end{equation}

\noindent Let us set $\delta_{n} = x_{n}-\beta z_{n}$, $\overline{\delta_{n}}=x_{n}+\beta z_{n}$, note  that $ Q_{\alpha}(u_{n}) = \delta_{n} \overline{\delta_{n}}$, so that $C_{n}^{\pm} = Q_{\alpha}(u_{n}) \pm \varepsilon$  and $B_{n}= -2(\beta \delta_{n} -x_{n}\e_{2n}(\beta)) $. The Dirichlet lattice point $u_{n}=(x_{n},0,z_{n})$ is exterior to $\{ -\varepsilon \leq Q_{\alpha}  \leq \varepsilon \}$, changing $u_n$ to $-u_{n}$ if necessary we can assume that $Q_{\alpha}(u_{n}) >\varepsilon$.  Thus, $x_{n}^{2} > \beta^{2}z_{n}^{2} \pm \varepsilon$ and  in particular $C_{n}^{\pm} > 0$.

\noindent We deduce from (\ref{dirich2}) and  (\ref{boundszn}) the following bounds for $ z_{n} $,$ \delta_{n}$ and $ x_{n} $, 
\begin{equation}\label{zn}
(C\q_{2n}^{1-\eta})^{1/\theta} = (CN_{n})^{1/\theta} < z_{n} \leqslant N_{n} = \q_{2n}^{1- \eta}.
\end{equation}

\begin{equation}\label{deltan}
\dfrac{C}{\q_{2n}^{(1-\eta)\theta}} = \dfrac{C}{N_{n}^{\theta}}   \leq  \dfrac{C}{z_{n}^{\theta}}< \delta_{n} \leq \dfrac{1}{N_{n}}= \dfrac{1}{\q_{2n}^{1-\eta}}.
\end{equation}

\noindent Since $x_{n} = \beta z_{n} + \delta_{n}$, using (\ref{zn}) and (\ref{deltan}) we get

\begin{equation}\label{xn}
\beta C^{1/\theta} \q_{2n}^{(1- \eta)/\theta} + \dfrac{C}{\q_{2n}^{(1- \eta)\theta}} <  \beta z_{n} + \dfrac{C}{N_{n}^{\theta}} < x_{n} \leq \beta z_n +\dfrac{1}{N_{n}} \leq \beta \q_{2n}^{1- \eta} + \dfrac{1}{\q_{2n}^{1- \eta}}.
\end{equation}

\noindent We define the following quantities, 

\begin{center}
$ U_{n} := -\dfrac{B_{n}}{A_{n}}$ and $ V_{n} := \dfrac{A_{n}}{B_{n}^{2}} $
\end{center}
thus $\Delta_{n}^{\pm}(\varepsilon)= B_{n}^{2} \left(1- 4V_{n}(Q_{\alpha}(u_{n}) \pm \varepsilon)\right).$  We are going to show that $V_{n}$ tends to zero, this will prove that $\Delta_{n}^{\pm}(\varepsilon)$ are both positive for $n$ large enough.

\noindent 

\noindent For any positive $n$, the inequalities in (\ref{error2}) gives 
 \begin{equation}\label{error5}
\dfrac{1}{2\q_{2n+1}(\beta)^{2}}  \leq    \e_{2n}(\beta) < \dfrac{1}{\q_{2n}(\beta)^{2}}.
\end{equation}

\noindent Concerning $B_{n}= 2\{ x_{n}\e_{2n}(\beta)- \beta \delta_{n}  \}$, with (\ref{deltan}), (\ref{xn}) and (\ref{error5}) one has that

\begin{equation}
\dfrac{1}{2\q_{2n+1}^{2}} \left( \dfrac{C}{N_{n}^{\theta}}  + \beta (CN_{n})^{1/\theta} \right)- \dfrac{2\beta}{N_{n}} < x_{n}\e_{2n}-2 \delta_{n} \beta \leq \dfrac{1}{\q_{2n}^{2}} \left( \beta  N_{n}+ \dfrac{1}{N_{n}} \right)- \dfrac{2\beta}{N_{n}^{\theta}}.
\end{equation}

\noindent We can rearrange the terms in order to get the following bounds for $B_n$

\begin{equation}\label{boundsBn}
\dfrac{1}{\q_{2n}} \left( \dfrac{C}{\lambda_{2n}^{2} \q_{2n}N_{n}^{\theta}} + \dfrac{ \beta (CN_{n})^{1/\theta}}{\lambda_{2n}^{2} \q_{2n}} - \dfrac{2\beta \q_{2n}}{N_{n}}\right) \leq B_n  \leq  \dfrac{2}{\q_{2n}^{2}} \left( \beta  N_{n}+ \dfrac{1}{N_{n}} \right)- \dfrac{4\beta}{N_{n}^{\theta}}.
\end{equation}

\noindent Replacing $N_{n}$ by $  \q_{2n}^{1-\eta} $ in (\ref{boundsBn}) we obtain a lower bound for $B_{n}$

\begin{equation}\label{lowerBn2}
\dfrac{1}{\q_{2n}} \left( \dfrac{C}{\lambda_{2n}^{2} \q_{2n}^{1+(1-\eta)\theta}} +  \beta C^{1/\theta}    \dfrac{\q_{2n}^{(1-\eta)/\theta -1}}{\lambda_{2n}^{2}}     -2\beta \q_{2n}^{\eta}\right) \leq B_n.
\end{equation}

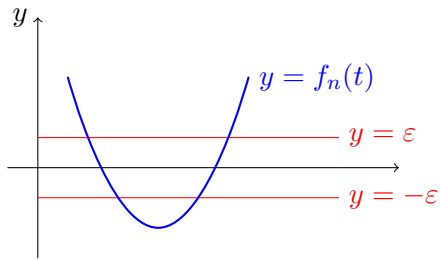
\begin{figure}

\begin{tikzpicture}[scale=0.4,cap=round]
 \tikzset{axes/.style={}}
 % The graphic
 \begin{scope}[style=axes]
 \draw[->,name path=c4] (-1,0) -- (12,0) node[below] {};
 \draw[->] (0,-3)-- (0,5) node[left] {$y$};
 \draw[red,name path=c1] (0,1)-- (10,1) node[right] {$y=\varepsilon$};
 \draw[red,-,name path=c2] (0,-1)-- (10,-1) node[right] {$y=-\varepsilon$};

 \draw [blue,thick,-,name path=c3] plot [smooth,tension=1] coordinates { 
   (1,3) (4,-2)  (7,3)}  node[right] {$y=f_{n}(t)$}  ;

\end{scope}
\end{tikzpicture}

\caption{ The domain  $-\varepsilon \leq  Q_{\alpha}(v_{n}(t))=f_{n}(t) \leq \varepsilon$ is supported by two intervals.}
\end{figure}

\noindent We claim that  $V_n$ tends to zero  as $n$ goes to infinity. Indeed, one has
$$V_{n} =\dfrac{\e_{2n}^{2}-2\beta \e_{2n}+ 1/\q_{2n}^{2}}{\left( x_{n}\e_{2n}-2\beta \delta_{n} \right)^{2}}.$$

\noindent For the numerator of $V_{n}$ we have the bound 
\begin{equation}
A_{n}= \e_{2n}^{2}-2\beta \e_{2n}+ \q_{2n}^{-2} \leq \dfrac{1}{\q_{2n}^{4}}- \dfrac{\beta}{\q_{2n+1}^{2}}+\dfrac{1}{\q_{2n}^{2}}= \dfrac{1}{\q_{2n}^{2}} \left( 1-  \dfrac{\beta}{\lambda_{2n}^{2}}+ \dfrac{1}{\q_{2n}^{2}}\right).
\end{equation}

\begin{figure}\label{figure2}

\def\tikzscale{0.4}
\begin{tikzpicture}[scale=\tikzscale]

\tikzset{
    elli/.style args={#1:#2and#3}{
        draw,
        shape=ellipse,
        rotate=#1,
        minimum width=2*#2,
        minimum height=2*#3,
        outer sep=0pt,
    }
}

%
% #1 optional parameters for \draw
% #2 angle of rotation in degrees
% #3 offset of center as (pointx, pointy) or (name-o-coordinate)
% #4 length of plus (semi)axis, that is axis which hyperbola crosses
% #5 length of minus (semi)axis
% #6 how much of hyperbola to draw in degrees, with 90 you‚Äö√Ñ√¥d reach infinity
%
\newcommand\tikzhyperbola[6][thick]{%
    \draw [#1, rotate around={#2: (0, 0)}, shift=#3]
        plot [variable = \t, samples=1000, domain=-#6:#6] ({#4 / cos( \t )}, {#5 * tan( \t )});
    \draw [#1, rotate around={#2: (0, 0)}, shift=#3]
        plot [variable = \t, samples=1000, domain=-#6:#6] ({-#4 / cos( \t )}, {#5 * tan( \t )});
}

\def\angle{90}
\def\bigaxis{3.2cm}
\def\smallaxis{1.5cm}

\draw [->, color=black, line width = 0.4pt, smooth] (-10, 0) -- (10, 0) node [right] {$x$};
\draw [->,color=black, line width = 0.4pt, smooth] (0, -10) -- (0, 16) node [above] {$z$};
\draw [dotted,color=black, line width = 1pt, smooth] (-9, 9) -- (16, 9) node [right] {$z= (CN_{n})^{1/\theta}= (C\q_{2n}^{1-\eta})^{1/\theta}$};

\draw [dotted,color=black, line width = 1pt, smooth] (-9, 11) -- (16, 11) node [right] {$z= N_{n}= \q_{2n}^{1-\eta}$};

\draw[domain=-7:7,scale=1,red] plot (\x,{2*\x -1.5});
\draw  (9,13) node[above] {$z=\beta^{-1}x$};
\draw  (12.5,11) node[above] {$ \color{red}(\mathcal{L}^{n}_{\beta}) : u_n + \mathbb{R}(c_{2n}, \dfrac{1}{\q_{2n}},1)$};
\draw[dotted] (5.8,9.3) node[above] {$\bullet$};

\draw[dotted] (-5,12.3) node[right] {$ \color{blue} Q_{\alpha} = - \varepsilon$};
\draw[dotted] (-6.5,12.3) node[left] {$  \color{blue} Q_{\alpha} = \varepsilon$};

\draw  (6.2,10) node[right] {$u_{n}=(x_{n}, 0, z_{n})$};
\coordinate (center) at (0, 0);

\node(center) (e) {};

%\draw [-{stealth}, line width = 1.2pt, color = orange] ([shift={(\angle:-12)}] e.center) -- ([shift={(\angle:12)}] e.center) node [above right] {$\bm{a}_1$};
%\draw [-{stealth}, line width = 1.2pt, color = orange] ([shift={(90+\angle:-8)}] e.center) -- ([shift={(90+\angle:8)}] e.center) node [above left]  {$\bm{a}_2$};

\tikzhyperbola[line width = 1.2pt, color=blue!80!black]{\angle}{(center)}{\bigaxis}{\smallaxis}{77}

\pgfmathsetmacro\axisratio{\smallaxis / \bigaxis}

% asymptotes
\def\lengthofasymptote{15}
\draw [color=black, line width = 0.4pt, rotate around={\angle + atan( \axisratio ): (center)}]
    ($ (-\lengthofasymptote, 0) + (center) $) -- ++(2*\lengthofasymptote, 0) ;
\draw [color=black, line width = 0.4pt, rotate around={\angle - atan( \axisratio ): (center)}]
    ($ (-\lengthofasymptote, 0) + (center) $) -- ++(2*\lengthofasymptote, 0) ;

\tikzhyperbola[line width = 1.2pt, color=blue!80!black]{90+\angle}{(center)}{\smallaxis}{\bigaxis}{76}

\end{tikzpicture}

\caption{In blue the level sets $Q_{\alpha}=\pm \varepsilon$ and in grey the generatrix $x=\beta z$ of the cone $Q_{\alpha}=0$ projected on the $xz$-plane. The red line $(\mathcal{L}^{n}_{\beta}) $ cuts $\{Q_{\alpha}=\pm \varepsilon\}$ in 4 points as $n$ gets large. The dotted lines represents the bounds for $z_{n}$.}
\end{figure}
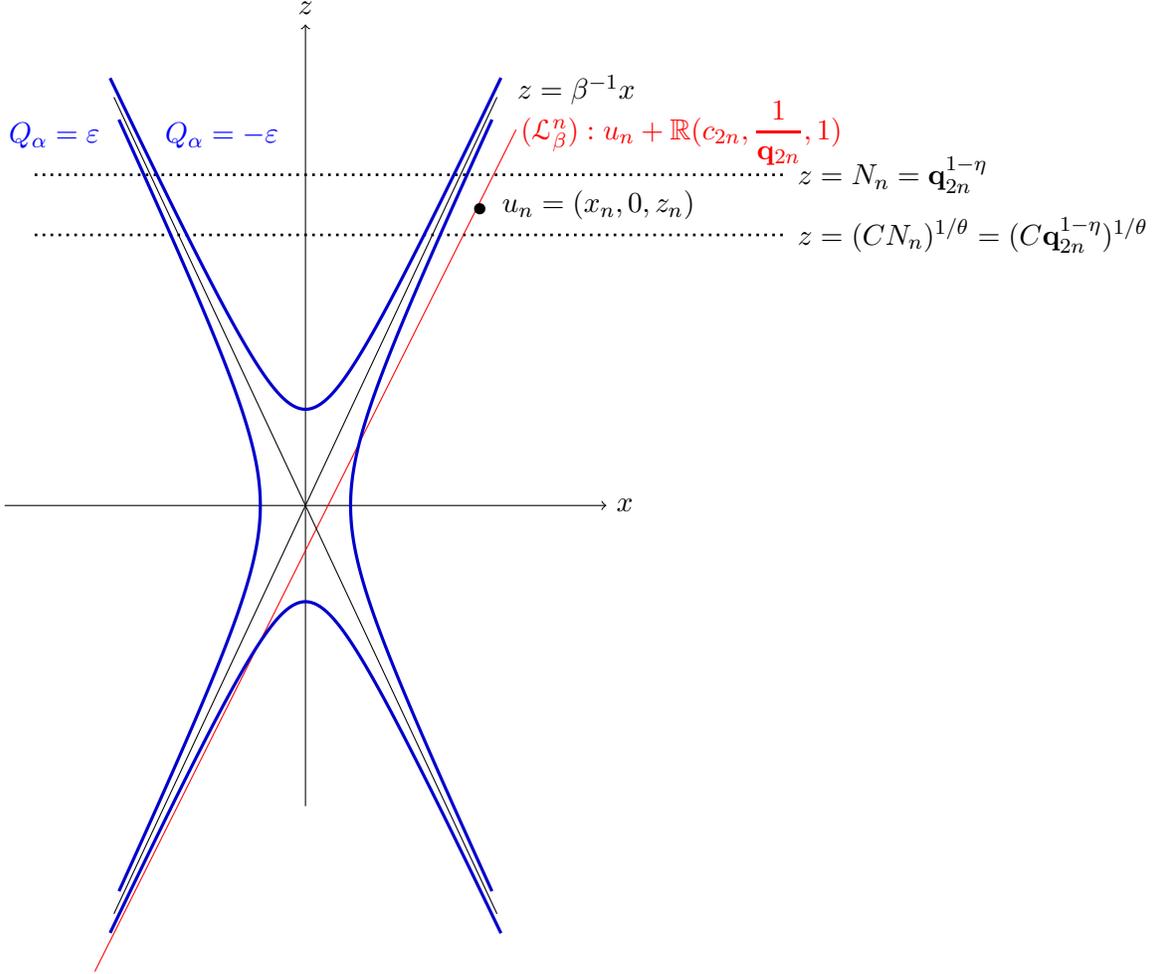

\noindent Thus, using (\ref{lowerBn2}) we get
\begin{equation}
0 < |V_n| \ll \dfrac{\left( 1-  \dfrac{\beta}{\lambda_{2n}^{2}}+ \dfrac{1}{\q_{2n}^{2}}\right)}{ \left( \dfrac{C}{\lambda_{2n}^{2} \q_{2n}^{1+(1-\eta)\theta}} +  \beta C^{1/\theta}  \dfrac{\q_{2n}^{(1-\eta)/\theta -1}}{\lambda_{2n}^{2}}     -2\beta \q_{2n}^{\eta} \right)^{2}}.
\end{equation}

\noindent Taking under consideration the fact that $ \lambda_{2n}^{2} \asymp \q_{2n}^{2(\mu -2)} $ which follows from (\ref{lambdamu}) we obtain that $$\lim_n V_n =0.$$

\noindent Using (\ref{ineqbadfinal}) we infer that 
\begin{equation}\label{Qn}
\vert Q_{\alpha}(u_{n}) \vert \leq 2C\beta + \dfrac{1}{N_{n}}=2C\beta + \dfrac{1}{\q_{2n}^{1-\eta}}.
\end{equation}

\noindent Thus the term $$1-4V_{n}(Q_{\alpha}(u_{n}) \pm \varepsilon)$$ can be made positive and less than 1  provided $n$ is taken large enough. Hence  the discriminants are always positive when $ n $ becomes larger than some positive integer $ n_{0} =n_{0}(\varepsilon)$ depending on $\varepsilon$. In this range the roots are given by
$$ t_{1,2,3,4}(n,\varepsilon)=\dfrac{1}{2}U_{n} \left( 1 \pm \sqrt{ 1-4V_{n}(Q_{\alpha}(u_{n}) \pm \varepsilon)} \right). $$

\noindent In more details, these correspond to the hitting times, $t_{1}< t_{2} < t_{3}< t_{4}$ given by
$$( (\mathcal{L}_{\beta}^{n})^{+} \cap \{ Q_{\alpha} = -\varepsilon\}) \left\lbrace \begin{array}{ccc}
t_{1}(n) &  = & \dfrac{1}{2}U_{n} \left( 1 - \sqrt{ 1-4V_{n}(Q_{\alpha}(u_{n})+\varepsilon)} \right) \\ \\

t_{4}(n) &  = &  \dfrac{1}{2}U_{n} \left( 1  +\sqrt{ 1-4V_{n}(Q_{\alpha}(u_{n})+\varepsilon)} \right)\\

\end{array}\right. $$
\noindent and 
$$( (\mathcal{L}_{\beta}^{n})^{+} \cap \{ Q_{\alpha} = \varepsilon\})  \left\lbrace \begin{array}{ccc}

t_{2}(n) &  = &  \dfrac{1}{2}U_{n} \left( 1  -\sqrt{ 1-4V_{n}(Q_{\alpha}(u_{n})-\varepsilon)} \right)\\ \\

t_{3}(n) &  = &  \dfrac{1}{2}U_{n} \left( 1 + \sqrt{ 1-4V_{n}(Q_{\alpha}(u_{n})-\varepsilon)} \right).\\
\end{array}\right. $$

\noindent In particular the values of  $t$ for which $v_{n}(t) \in (\mathcal{L}_{\beta}^{n})^{+} \cap \mathcal{A}(\varepsilon)$ is the union the two disjoint intervals $I_{1}^{n}(\varepsilon)=[t_{1}, t_{2}]$ and $I_{2}^{n}(\varepsilon)=[t_{3}, t_{4}]$ when $n \geq n_{0}$.  The proposition \ref{exit times} is proved.

\begin{flushright}
$\square$
\end{flushright}

\noindent \textit{Remark.} The signs and the hitting times are not important for our purposes.  The two intervals of the proposition \ref{exit times},  $I_{1}^{n}$ and $I_{2}^{n}$ plays a symmetric role and their size is the same. One interval should comprise negative times while the other consists of positive ones. Assume for instance, changing the order if necessary, that the roots are sorted such that $t_1 <t_2  < t_3 < t_4$ for $n \geqslant n_{0}$. Then in view of the previous proposition, for $n \geqslant n_{0}$ the intersection of the half-line $(\mathcal{L}_{\beta}^{n})^{+}$ with the two level sets $\{ Q_{\alpha}=\pm \varepsilon \}$ behaves as follows, see figure \ref{figure2}
$$\left\lbrace \begin{array}{ccc}
 Q_{\alpha}(v_{n}(t)) < -\varepsilon  &  \mathrm{if} & 0 < t < t_{1} \ \ \ (\mathrm{out}) \\
-\varepsilon \leq Q_{\alpha}(v_{n}(t)) \leqslant \varepsilon &  \mathrm{if} & t_1 \leqslant t < t_{2} \ \ \ (\mathrm{in})\\
\varepsilon < Q_{\alpha}(v_{n}(t)) &  \mathrm{if} & t_{2} \leqslant t < t_{3} \ \ \ (\mathrm{out}) \\
 -\varepsilon < Q_{\alpha}(v_{n}(t)) \leq \varepsilon &  \mathrm{if} & t_3 \leqslant t < t_{2} \ \ \ (\mathrm{in})\\
 Q_{\alpha}(v_{n}(t)) <-\varepsilon &  \mathrm{if} & t_{4}< t \ \  \ \ \ \ (\mathrm{out}).
\end{array}\right. $$

\section{A Solution to the Oppenheim conjecture for $Q_{\alpha}$}

Let $\varepsilon>0$ be an arbitrary small real number, we are interested to finding $n$ and $t_{n}$ such that $v_n(t_{n})$ is a nonzero vector is in $\mathbb{Z}^{3} \cap \mathcal{A}(\varepsilon) $ i.e. 
$$0< \vert Q_{\alpha}(v_{n}(t_{n})) \vert \leq \varepsilon. $$
In order that  $v_n(t_{n})$ provides the required lattice point, we necessarily need  $t_{n}$ to be a multiple of $\q_{2n}$ so that we can clear the denominators. By symmetry, we only need to focus on one interval, say $I_{1}^{n}(\varepsilon)$.  The following combinatorial argument shows that it is always possible to do so for large enough values of $n$.

\begin{lem} \label{multiple} 
There exists a positive  integer $n_{1}(\varepsilon)$ such that the interval $I_{1}^{n}(\varepsilon)=[t_{1}, t_{2}]$ contains a multiple of $\q_{2n}$ whenever $ n \geq n_{1}$.
\end{lem}

\noindent \textbf{Proof.} Let us set for each positive integer $ n $, the following counting function
\begin{center}
$  M_{n} := \mathrm{Card}  \left(   [t_{1}, t_{2}] \cap \mathbb{Z}\q_{2n} \right). $
\end{center}
 $  M_{n}$ is the number of multiples of $ \q_{2n} $ in $I_{1}^{n}=[t_{1}, t_{2}]$. In particular we have that $ M_{n} = \lfloor  \dfrac{l(I_{1}^{n})}{ \q_{2n}} \rfloor$ where $ l(I_{1}^{n}) $ is the length of the interval  $I_{1}^{n}$.
  The aim is to show that this quantity is $ \geq 1$ when $n$ is larger that a certain threshold $n_{1}$.  As $n$ gets large, we have that
 \begin{equation} \label{t1}
 t _{1} \asymp  2 U_{n} V_{n} (Q_{\alpha}(u_{n})-\varepsilon) = -\frac{2(Q_{\alpha}(u_{n})-\varepsilon) }{B_{n}} 
\end{equation}
and 
 \begin{equation} \label{t2}
 t _{2} \asymp  2 U_{n} V_{n} (Q_{\alpha}(u_{n})+\varepsilon) = -\frac{2(Q_{\alpha}(u_{n})+\varepsilon) }{B_{n}}.
\end{equation}

\noindent The length of the interval $I_{n}^{1}$ is asymptotically given by 

$$ l(I_{1}^{n}) =   |t_{2} - t_{1}|  \asymp  \dfrac{4\varepsilon}{\vert B_{n} \vert }.$$

\noindent Thus, 
$$ M_n \asymp \dfrac{4\varepsilon}{\q_{2n} |B_n|}.  $$

\noindent Concerning the denominator, 

$$\q_{2n} \vert B_n \vert  \leq \dfrac{1}{\q_{2n}}\left( \beta N_{n} + \dfrac{1}{N_{n}}\right) -\dfrac{2\beta}{N_{n}^{\theta}}\q_{2n}$$

$$\ \ \ \ \ \leq \beta \dfrac{N_{n}}{\q_{2n}} + \dfrac{1}{\q_{2n}N_{n}} -2\beta \dfrac{\q_{2n}}{N_{n}^{\theta}}$$

$$\ \ \ \ \ \  \leq  \dfrac{\beta}{\q_{2n}^{\eta}} + \dfrac{1}{\q_{2n}^{2-\eta}} - \dfrac{2\beta}{\q_{2n}^{(1-\eta)\theta -1}}.$$

\noindent Then the choice\footnote{This is the only moment we need that $  \theta>1$} of $\eta$ gives that $(1-\eta)\theta -1 = \eta$, then 

$$\q_{2n} \vert B_n  | \leq  \dfrac{\beta}{\q_{2n}^{\eta}} + \dfrac{1}{\q_{2n}^{2-\eta}} - \dfrac{2\beta}{\q_{2n}^{\eta}}.$$

$$\ \ \ \ \ \ \ \  \ \ \ \ \ \ \   \leq  \dfrac{\beta}{\q_{2n}^{\eta}} \left| 1 - \dfrac{1}{\beta \q_{2n}^{2(1-\eta)}}  \right|.$$

\noindent Thus we have the upper estimate 

$$ \q_{2n} |B_{n} | \ll \dfrac{1}{\q_{2n}^{\eta}}.$$

\noindent Taking the inverse, 

\begin{equation} \label{mn1}
\varepsilon \q_{2n}^{\eta} \ll \dfrac{4\varepsilon}{\q_{2n} |B_n|}.
\end{equation} 

\noindent Finally we infer the following crucial bound

\begin{equation} \label{mn2}
\varepsilon \q_{2n}^{\eta} \ll M_n
\end{equation}

\noindent  In particular, since $(\q_{2n}^{\eta})_{n}$ diverges there exists $n_{1}(\varepsilon)$ such that for all $n \geq n_{1}(\varepsilon)$

\begin{equation}
1<\varepsilon\q_{2n}^{\eta}.
\end{equation}

Hence (\ref{mn2}) shows that  $M_{n} \geq 1$ for $n \geq n_1$, meaning that the interval of times $I_{1}^{n}(\varepsilon)=[t_{1}, t_{2}]$ contains at least one multiple of $\q_{2n}$ for $n \geq n_1$. Let us estimate the integer $n_{1}$ which depends on the choice of $ \varepsilon $ and $ \eta $ , and which can be seen formally as
$$\varphi_{\beta}( \varepsilon) =  \min\{ n\geq n_{0} \ \ \vert \ \  M_{n} \geq 2  \}.  $$

Here $n_{0}=n_{0}(\varepsilon)$ is the least integer which ensures that $I_{1}^{n}(\varepsilon) \neq \emptyset$ coming from Proposition  \ref{exit times} whereas $n_{1}(\varepsilon)$ is the least integer such that $I_{1}^{n}(\varepsilon)$ contains a multiple of $\q_{2n}$. In particular, $n_{0}(\varepsilon) < n_{1}(\varepsilon)$. 
The number $ n_{1} $ is not going to be optimal, i.e. it will be an upper estimate for $ \varphi_{\beta}(\varepsilon) $.

\noindent Since $2^{n-1} \leq \q_{2n}  $, a sufficient condition in order the inequality  $  \varepsilon \q_{2n}^{\eta} > 1$ to hold is

$$  \varepsilon  2^{\eta(n-1)} > 1.$$

\noindent Applying logarithms, we get

$$ n > 1+ \eta^{-1}  \left| \ln \left( \dfrac{1}{\varepsilon} \right) \right| /\ln 2 .$$

\noindent Thus a good choice for $n_{1}$ is

$$ n_{1}(\varepsilon) := 2+  \lfloor \eta^{-1}  \left| \ln \left(\varepsilon \right) \right| /\ln 2 \rfloor.$$
\noindent This finishes the proof of the Lemma. 
\begin{flushright}
$\square$
\end{flushright}

\subsection*{Proof of  Theorem \ref{main}.}
Let $ \varepsilon >0 $ be fixed. 

\noindent \textbf{Case 1} Assume $\beta$ is a Liouville number and let $n$ be a positive large enough integer so that  $$2^{-n} \beta+ 2^{-2(n+2)} \leq \varepsilon.$$ Since $\mu(\beta)=\infty$ and given $n$ as above we can always find a rational number $p/q$ such that 
\begin{equation}\label{trans1}
 \left|q \beta - p \right| < \dfrac{1}{q^{n+2}}.
\end{equation}

\noindent From this, we deduce that $ q\beta - 1/q^{n+2} < p < q\beta - 1/q^{n+2} $, thus 

\begin{equation}\label{trans2}
2q\beta - 1/q^{n+2} < p  + \beta  q < 2q\beta +  1/q^{n+2}.
\end{equation}

\noindent Thus, 

\begin{equation}
|Q_{\alpha}(p,0,q)| = |p^2 - \beta^{2} q^{2}|= (p + \beta q) |p - \beta q|< \dfrac{1}{q^{n+2}} (2q\beta +  1/q^{n+2}).
\end{equation}

\noindent Since $q \geq 2$, 
$$ |Q_{\alpha}(p,0,q)| < \dfrac{1}{2^{n+2}} (4\beta +  1/2^{n+2}) .$$

\noindent The choice of $n$ implies that $v=(p,0,q)$ is a nonzero integral solution of 
 $$ |Q_{\alpha}(p.0.q)| < \varepsilon.$$ 
\noindent In other words, the Oppenheim conjecture holds for $Q_{\alpha}$ in this case.

\noindent \textbf{Case 2} Assume $\beta$ is a not a Liouville number,

The lemma (\ref{multiple}) shows that there exists an explicit  integer $ n_{1}(\varepsilon) >0$ such that  $ I_{1}^{n}(\varepsilon)$ contains a multiple of $\q_{2n_{1}} $ say $a_{n_{1}} \q_{2n_{1}}  \in [t_{1}, t_{2}]$ where $a_{n_{1}}$ is a nonzero integer. The proposition (\ref{exit times}) implies that $  v_{n_{1}}(a_{n_{1}}\q_{2n_{1}}) \in \mathcal{A}(\varepsilon)$. Moreover, 
$$v_{n_{1}}(a_{n_{1}}\q_{2n_{1}}) = \left( x_{n_{1}}-a_{n_{1}}p_{2n_{1}},  -a_{n_{1}},  z_{n_{1}}- a_{n_{1}} \q_{2n_{1}} \right) \in \mathbb{Z}. $$
Thus, we have a nonzero integral vector $v_{1} := v_{n_{1}}(a_{n_{1}} \q_{2n_{1}})$ in $   \mathcal{A}(\varepsilon)$, that is, 
$$ \vert Q_{\alpha}\left( v_{1}   \right)\vert  \leqslant \varepsilon. $$
This proves the first assertion of the theorem. We give an estimate the size of the solution, set
$$ \| v_1 \|_{\infty} = \max \{ \vert x_{n_{1}}-a_{n_{1}}p_{2n_{1}} \vert, \vert z_{n_{1}}-a_{n_{1}}p_{2n_{1}} \vert, \vert a_{n_{1}} \vert  \}.$$

\noindent A crude bound is given by
$$ \| v_1 \|_{\infty}  \leq  \vert x_{n_{1}} \vert + \vert a_{n_{1}} \vert   p_{2n_{1}}.$$
We know from (\ref{xn}) that 
$$ |x_{n_{1}}| \lesssim \beta \q_{2n_{1}}^{1-\eta}. $$

\noindent Also by construction we have $a_{n_{1}} \in [\dfrac{t_{1}}{\q_{2n_{1}}} ,  \dfrac{t_{2}}{\q_{2n_{1}}} ]$, thus in view of (\ref{Qn}), (\ref{t1}) and (\ref{t2}) one has 
$$ \vert  a_{n_{1}} \vert \leq  \dfrac{|Q_{\alpha}(u_{n})| + \varepsilon}{\q_{2n_{1}} |B_{n_{1}}|}  \leq \dfrac{ 2C \beta  + \q_{2n}^{-1+\eta} + \varepsilon}{\q_{2n_{1}} |B_{2n_{1}}|}. $$
Thus, 
$$  \| v_1 \|_{\infty} \lesssim \q_{2n_{1}}^{1-\eta}+ \dfrac{ 2C \beta  + \varepsilon}{\q_{2n_{1}} |B_{n_{1}}|}  \mathbf{p}_{2n_{1}} + \dfrac{c_{2n_{1}}}{\q_{2n_{1}}^{1-\eta} |B_{n_{1}}|}$$

\noindent or equivalently

$$  \| v_1 \|_{\infty} \lesssim \q_{2n_{1}}^{1-\eta}+ \dfrac{ 2C \beta  + \varepsilon}{|B_{n_{1}}|} \mathbf{c}_{2n_{1}}+  \dfrac{c_{2n_{1}}}{\q_{2n_{1}}^{1-\eta} |B_{n_{1}}|}.$$

\noindent Using (\ref{lowerBn2})
\begin{equation}\label{bn1}
\dfrac{1}{  \mid  B_{n_{1}}  \vert } \leq  \dfrac{1}{ \vert C{\lambda_{2n_{1}}^{-2} \q_{2n_{1}}^{-2-(1-\eta)\theta}+  \beta C^{1/\theta} \lambda_{2n_{1}}^{-2} \q_{2n_{1}}^{(1-\eta)/\theta -2} -2\beta \q_{2n_{1}}^{\eta-1}\mid }}.
\end{equation}

\noindent We have the relation $  1-\eta =2/(\theta(\theta+1))$ and the irrationality measure $\mu$ which comes into play using (\ref{lambdamu}), thus

$$\lambda_{2n}^{-2}\q_{2n}^{-2-(1-\eta)\theta} \asymp \q_{2n}^{-2(\mu-2)-2-(1-\eta)\theta} = \q_{2n}^{-2\mu +2 -2/ (\theta+1)}=  \q_{2n}^{-2(\mu -1 +1/ (\theta+1))}$$

and 
$$\lambda_{2n}^{-2}\q_{2n}^{(1-\eta)/\theta-2} \asymp \q_{2n}^{-2(\mu-2)+(1-\eta)/\theta-2} = \q_{2n}^{-2\mu +2 + 2/(\theta^{2}(\theta+1))}= \q_{2n}^{-2(\mu -1 - 1/(\theta^{2}(\theta+1)) } .$$

\noindent Therefore, 

$$\dfrac{1}{  \mid  B_{n_{1}}  \vert } \ll  \dfrac{1}{ \vert C\q_{2n}^{-2(\mu -1 +1/ (\theta+1))} +  \beta C^{1/\theta}\q_{2n}^{-2(\mu -1 - 1/(\theta^{2}(\theta+1)) } -2\beta \q_{2n_{1}}^{-2/(\theta(\theta+1))}\mid }.$$\\

\noindent Let us set   $\kappa = \dfrac{2}{\theta(\theta +1)}= 1-\eta $,

$$\dfrac{1}{  \mid  B_{n_{1}}  \vert } \ll  \dfrac{1}{ \vert C{\q_{2n_{1}}^{-2(\mu -1)-\kappa}  +  \beta C^{1/\theta}\q_{2n_{1}}^{-2(\mu-1)+\kappa/\theta} -2\beta \q_{2n_{1}}^{-\kappa}\mid }}.$$\\

$$ \ \ \ \ \  \ \ \ \  \ll  \dfrac{\q_{2n_{1}}^{\kappa}}{ \vert 1- C {\q_{2n_{1}}^{-2(\mu -1)}/2\beta  - C^{1/\theta} \q_{2n_{1}}^{-2(\mu-1)+\kappa/\theta +\kappa}/2  \mid }}.$$\\

\noindent Thus,  since $2(\mu-1)\geq 2>\kappa/\theta +\kappa$, one has

$$\dfrac{1}{  \vert  B_{n_{1}}  \vert }  \ll \q_{2n_{1}}^{\kappa}= \q_{2n_{1}}^{1-\eta} .$$

\noindent Hence we get,

$$  \| v_1 \|_{\infty} \ll  \q_{2n_{1}}^{1-\eta} + O(1).$$

\noindent In short, 
$$  \| v_1 \|_{\infty} \ll \q_{2n_{1}}^{1-\eta}.$$

\noindent By definition we have  $\eta = 1-2/(\theta +1) $ so that the last inequality reads

$$  \| v_1 \|_{\infty} \ll  \q_{2n_{1}}^{2/(\theta + 1)}$$

\noindent where 
$$  n_{1}(\varepsilon) := 2+  \lfloor \eta^{-1}  \left| \ln \left(\varepsilon \right) \right| /\ln 2 \rfloor.$$

\noindent This finishes the proof of Theorem \ref{main}.
\begin{flushright}
$ \square$
\end{flushright}

\subsection*{Proof of  Corollary \ref{cor}.}
$(1)$ By assumption there exists $\gamma \in \mathrm{SL}(3, \mathbb{Q})$, such that $Q(x)= Q_{\alpha}(\gamma x)$. Let us consider an arbitrary real $\varepsilon > 0$.  Let $a= \mathrm{lcm} \{ \mathrm{den}((\gamma^{-1})_{i,j}), 1\leqslant i,j \leqslant 3\} $ be the least common multiple of the denominator of the coefficients of $\gamma^{-1}$. Thus $a\gamma^{-1}$ is an integral matrix.
Theorem \ref{main}  gives the existence of a nonzero integral vector $v \in \mathbb{Z}^{3}$ such that $|Q_{\alpha}(v)| \leq \varepsilon/a^{2}$. Since $a\gamma^{-1} \in \mathrm{SL}(3, \mathbb{Z})$,  $v_{1}= a\gamma^{-1}v$ is a nonzero integral vector such that 
$$ |Q(v_{1})|= |Q(a\gamma^{-1}v)|= a^{2}|Q(\gamma^{-1}v)  |= a^{2}|Q_{\alpha}(v)| \leq \varepsilon.$$
Hence $Q$ satisfies the Oppenheim conjecture. 

\noindent $(2)$ Let $h \in H$ such that $Q(x)=Q_{0}(hx)$, where for some $A \in \mathrm{SL}(3,\mathbb{Q}) $ and $ h_{33} \notin \mathbb{Q}$ one has $$  h =\left[\begin{array}{c|c}A & 0 \\\hline 0 & h_{33}\end{array}\right]. $$

The matrix $h$ can factorized as follows
$$h=  \left[\begin{array}{c|c} I_{2} & 0 \\\hline 0 & h_{33}\end{array}\right]   \left[\begin{array}{c|c} A & 0 \\\hline 0 &1\end{array}\right].$$

Set $\alpha=h_{33}^{2}$ and $\gamma =  \left[\begin{array}{c|c} A & 0 \\\hline 0 &1\end{array}\right] \in \mathrm{SL}(3,\mathbb{Q})$, thus
$$ Q(x)=Q_{0}(hx) = Q_{\alpha}(\gamma x).$$

The form $Q$ is $\mathrm{SL}(3,\mathbb{Q}) $-equivalent to the form $Q_{\alpha}$. Then the assertion $(1)$ of the corollary allows us to show that $Q$ fullfills the conjecture. 
\begin{flushright}
$\square$
\end{flushright}

\end{document}